\documentclass{amsart}
\usepackage{amssymb}
\usepackage{amsfonts}
\usepackage{latexsym}

\newtheorem{theorem}{Theorem}[section]
\newtheorem{lemma}[theorem]{Lemma}
\newtheorem{proposition}[theorem]{Proposition}
\newtheorem{corollary}[theorem]{Corollary}
\theoremstyle{definition}
\newtheorem{definition}[theorem]{Definition}
\newtheorem{example}[theorem]{Example}
\newtheorem{conjecture}[theorem]{Conjecture}
\newtheorem{remark}[theorem]{Remark}

\newcommand{\Rep}{\text{Rep}}

\newcommand{\raro}{\to}

\newcommand{\ot}{\otimes}
\newcommand{\qexp}{{\rm {qexp}}}
\newcommand{\ben}{\begin{enumerate}}
\newcommand{\een}{\end{enumerate}}

\newcommand{\Rad}{{\rm {Rad}}}

\begin{document}

\title{On the quasi-exponent of finite-dimensional Hopf algebras}

\author{Pavel Etingof}
\address{Department of Mathematics, Massachusetts Institute of Technology,
Cambridge, MA 02139, USA}
\email{etingof@math.mit.edu}

\author{Shlomo Gelaki}
\address{Department of Mathematics, Technion-Israel Institute of
Technology, Haifa 32000, Israel}
\email{gelaki@math.technion.ac.il}


\maketitle


\begin{section}
{Introduction}
In [EG] we introduced and studied
a new invariant of Hopf algebras -- the exponent.
This invariant generalizes the exponent of a group (the least common
multiple of orders of its elements).
It has a number of interesting properties, and has found some
applications [KSZ].

More specifically, consider Hopf algebras over $\mathbb{C}.$
Then the exponent has a very nice behavior for
{\it semisimple} Hopf algebras. Namely, it is shown in [EG]
that in this case the exponent is finite and divides
the cube of the dimension. On the contrary,
for a non-semisimple Hopf algebra, the exponent
is usually (and possibly always) infinite.

Thus, one is tempted to look for a more refined invariant
of Hopf algebras, which would coincide with the
exponent in the semisimple case, but would be finite
if the Hopf algebra is finite-dimensional. This problem is
solved in the present paper, by introduction
of the {\it quasi-exponent}.

Recall from [EG] that the exponent of
a finite-dimensional Hopf algebra $H$
is the order of the Drinfeld element $u$
of the Drinfeld double $D(H)$ of $H.$
Recall also that while
this order may be infinite,
the eigenvalues of $u$ are always roots of unity
([EG, Theorem 4.8]); i.e., some power
of $u$ is always unipotent. We are thus naturally led
to define the {\it quasi-exponent} of a finite-dimensional Hopf algebra
$H$ to be the order of unipotency of $u$.

The goal of the paper is to create a theory of quasi-exponent, which
would be parallel to the theory of the exponent developed in [EG]. In
particular, similarly to [EG], we give two other equivalent
definitions of the quasi-exponent, and prove that it is invariant
under twisting. Furthermore, we prove that the quasi-exponent of a
finite-dimensional pointed Hopf algebra $H$ is equal to the exponent of
the group $G(H)$ of grouplike elements of $H.$ (In particular,
the order of the squared antipode of $H$ divides $\exp(G(H)).$)
As an application, we find that if $H$ is obtained
by twisting the quantum group at root of unity $U_q({\mathfrak g})$
then the order of any grouplike element in $H$ divides the order
of $q$.

In this paper, all Hopf algebras will be over the field of
complex numbers $\mathbb{C},$
unless otherwise is specified.

\end{section}

\begin{section}
{Definition and elementary properties of quasi-exponent}

Let $H$ be a complex finite-dimensional Hopf algebra with multiplication map $m,$
comultiplication map $\Delta,$ unit map $\iota,$ counit map $\varepsilon$ and antipode
$S.$ Recall that the Drinfeld double
$D(H)=H^{*cop}\ot
H$ of $H$ is a quasitriangular Hopf algebra with universal $R-$matrix $R=\sum_ih_i\ot
h_i^*,$ where $\{h_i\},$ $\{h_i^*\}$ are dual bases for
$H$ and $H^*$ respectively. Let
\begin{equation}\label{dre}
u:=m_{21}(Id\ot S)(R)=\sum_iS(h_i^*)h_i,
\end{equation}
be the Drinfeld element of
$D(H),$ where $Id$ denotes the identity map $H\to H,$
$S$ is the antipode of $D(H)$ and $m_{21}$ is the composition of the usual
permutation map with the multiplication map $m$ of $D(H).$ Sometimes we shall use the
notation $u_H$ for emphasis. By [D],
\begin{equation}\label{ssq}
S^2(x)=uxu^{-1},\;x\in D(H).
\end{equation}

We start by making the following definition.

\begin{definition}\label{qexp}
Let $H$ be a finite-dimensional Hopf algebra over $\mathbb{C}.$
The quasi-exponent of $H,$ denoted by $\qexp(H),$ is the smallest positive
integer $n$ such
that $u^n$ is unipotent; that is, the smallest positive integer $n$ such that
$(1-u^n)^N=0$ for some positive integer $N.$
\end{definition}

\begin{remark}\label{pch}
Let $H$ be a finite-dimensional Hopf algebra over $\mathbb{C}.$ Then the following hold:
\ben
\item
If $H$ is semisimple, $\qexp(H)=\text{exp}(H)$ by [EG, Theorem 4.3].

\item
By [EG, Theorem 4.8], $\qexp(H)<\infty.$

\item
If $\text{exp}(H)<\infty$ then $\qexp(H)=\text{exp}(H).$
\een
\end{remark}

Let $m_0=\varepsilon,$ $m_1=Id,$ $\Delta_0=\iota$ and $\Delta_1=Id,$ and for any
integer $n\ge 2$ let $m_n:H^{\ot n}\raro H$ and
$\Delta_n:H\raro H^{\ot n}$ be defined by $m_n=m\circ(m_{n-1}\ot
Id),$ and $\Delta_n=(\Delta_{n-1}\ot Id)\circ\Delta$. Let $T_n:H\rightarrow H$ be the
linear map given by
$$T_n(h)=m_n\circ(Id\ot S^{-2}\ot\cdots\ot S^{-2n+2})\circ\Delta_n(h)$$
and set
$$R_n:=R(Id\ot S^2)(R)\cdots(Id\ot S^{2n-2})(R).$$
\begin{proposition}\label{eqdef}
Let $H$ be a finite-dimensional Hopf algebra, and $f(t)=\sum_{i=0}^ma_it^i$
be a polynomial in $\mathbb{C}[t].$ The following conditions are equivalent:

\ben
\item $f(u)=0.$

\item $\sum_{i=0}^ma_iT_i=0.$

\item $\sum_{i=0}^ma_iR_i=0.$
\een
\end{proposition}
\begin{proof}
$(3)\Rightarrow (1):$ By a simple induction on $i,$ $m_{21}(Id\ot S)(R_i)=u^i$ (here we
use equation (\ref{ssq})). Hence (1) follows by applying $m_{21}(Id\ot S)$ to the
equation $\sum_{i=0}^ma_iR_i=0.$

$(2) \Leftrightarrow (3):$ It is straightforward to show that the Hexagon
axiom for $R$ implies that $(\Delta_n\ot Id)(R)=R_{1,n+1}R_{2,n+1}\cdots R_{n,n+1}.$
Therefore, $(T_i\ot Id)(R)=R_i.$  Thus applying (2) to the first component of $R$ we get (3).
Conversely, the equation $(\sum_{i=0}^ma_iT_i\ot Id)(R)=0$ implies (2).

$(1)\Rightarrow (3):$  Since $m_{21}(Id\ot S)(R_i)=u^i,$ one has $m_{21}(Id\ot
S)(\sum_{i=0}^ma_iR_i)=0.$ But the multiplication map $H^*\ot H\to D(H)$ is a linear
isomorphism [D], so $\sum_{i=0}^ma_iR_i=0,$ as desired.
\end{proof}

Proposition \ref{eqdef} allows us to give two equivalent definitions of the
quasi-exponent.
\begin{corollary}\label{eqdef1}
Let $H$ be a finite-dimensional Hopf algebra over $\mathbb{C}.$ Then:
\ben
\item $\sum_{k=0}^N (-1)^k \left(\begin{array}{c} N \\ k \end{array} \right)
R_{nk}=0$ for some positive integer $N$ if and only if $n$ is divisible by $\qexp(H).$
\item $\sum_{k=0}^N (-1)^k \left(\begin{array}{c} N \\ k \end{array} \right)
T_{nk}=0$ for some positive integer $N$ if and only if $n$ is divisible by $\qexp(H).$
\een
Thus, $\qexp(H)$ may be defined as the smallest positive integer satisfying either (1)
or (2).
\end{corollary}

Let us now list some of the elementary properties of $\qexp(H).$

\begin{proposition}\label{elem}
Let $H$ be a finite-dimensional Hopf algebra over $\mathbb{C}.$ Then:
\ben
\item
The quasi-exponents of Hopf subalgebras and quotients of $H$ divide $\qexp(H).$
\item
The order of any grouplike element $g$ of $H$ divides $\qexp(H).$
\item
$\qexp(H^*)=\qexp(H).$
\item
$\qexp(H_1\ot H_2)$ equals the least common multiple of $\qexp(H_1)$
and $\qexp(H_2).$
\item
The order of $S^2$ divides $\qexp(H)$ (in particular, $u^{\qexp(H)}$ is central in
$D(H)$).
\item
If $\qexp(H)=1$ then $H=\mathbb{C}.$
\item
$\qexp(H^{*cop})=\qexp(H).$
\een
\end{proposition}
\begin{proof}
(1) Follows from part (2) of Corollary \ref{eqdef1}.

\noindent
(2) The order of $g$ divides $\exp(G(H))=\exp(\mathbb{C}[G(H)])$ which divides $\qexp(H)$
by part (1).

\noindent
(3) Clear from part (2) of Corollary \ref{eqdef1}.

\noindent
(4) Clearly, $u_{H_1\ot H_2}=u_{H_1}\ot u_{H_2}.$ Hence if $u_{H_1}^n,$ $u_{H_2}^n$ are
both unipotent then so is $u_{H_1\ot H_2}^n.$ Conversely, if $u_{H_1\ot H_2}^n$ is
unipotent, so
is the operator
$u_{H_1\ot H_2}^n$ restricted to $D(H_1)\ot \mathbb{C}\lambda,$ where $\lambda$ is a non-zero
left integral of $D(H_2).$
But this operator equals $u_{H_1}^n$ since $\varepsilon(u)=1.$ Similarly, $u_{H_2}^n$ is
unipotent as well.

\noindent
(5) Let $n:=\qexp(H).$ Since $S^{2n}=\text{Ad}u^n,$ $S^{2n}$ is unipotent. However, it is
semisimple by [R], as it
is of finite order. Thus, $S^{2n}=Id.$

\noindent
(6) By part (5), $S^2=Id,$ so by [LR], $H$ is semisimple. But then $1=\qexp(H)=\exp(H)$
which implies that $Id=\varepsilon,$ and hence that $H=\mathbb{C},$ as desired.

\noindent
(7) Since $(D(H^{*cop}),R)\cong (D(H)^{op},R_{21})$ as quasitriangular Hopf algebras,
it follows that $u_H=u_{H^{*cop}}.$
\end{proof}

\begin{example}
Let $H$ be Sweedler's $4-$dimensional non-semisimple Hopf algebra generated by the
grouplike element $g$ and the skew-primitive element $x$ with $\Delta(x)=x\ot g+1\ot x$
[Sw]. Using the basis $\{1,g,x,gx\}$ of $H$ it is straightforward to
verify that
$\varepsilon -2m_2(Id\ot S^{-2})\Delta_2 + m_4(Id\ot S^{-2}\ot S^{-4}\ot
S^{-6})\Delta_4=0.$ Therefore, $\qexp(H)$ divides $2.$
Since it
is not equal to $1,$ we have $\qexp(H)=2=\exp(G(H)).$ In Section 4 we will generalize this
to any finite-dimensional pointed Hopf algebra.
\end{example}

\end{section}

\begin{section}
{Invariance of quasi-exponent under twisting}
In this section we study the invariance of $\qexp(H)$ under twisting. But first let us
recall Drinfeld's notion of a twist for Hopf algebras for the convenience of the reader.

\begin{definition}\label{twist}
Let $H$ be a Hopf algebra. A twist for $H$
is an invertible element $J\in H\ot H$ which satisfies
$$
(\Delta\ot Id)(J)(J\ot 1)=(Id\ot \Delta)(J)(1\ot J)\;\; and \;\;
(\varepsilon\ot Id)(J)=(Id\ot \varepsilon)(J)=1.
$$
\end{definition}
Given a twist $J$ for $H,$ one can construct a new Hopf algebra $H^J,$
which is the same as $H$ as an algebra, with coproduct $\Delta^J$ given by
$$\Delta^J(x)=J^{-1}\Delta(x)J,\; x\in H,$$
and antipode $S^J$ given by
$$S^J(x)=Q^{-1}S(x)Q,\; x\in H,$$
where $Q=m\circ(S\ot Id)(J)$ and $Q^{-1}=m\circ(Id\ot S)(J^{-1}).$

If $(A,R)$ is quasitriangular then so is $A^J$ with the universal $R$-matrix
$$R^J:=J_{21}^{-1}RJ.$$
In particular, since $H$ is a Hopf subalgebra of $D(H),$ we can twist
$D(H)$ using the twist $J,$ considered as an element of $D(H)\ot D(H),$ and obtain
$(D(H)^J,R^J).$
In [EG, Proposition 3.2] we proved that there exists a canonical isomorphism of
quasitriangular Hopf algebras between $(D(H)^J,R^J)$ and $(D(H^J),R^{(J)})$ (where
$R^{(J)}$ is the universal $R-$matrix of $D(H^J)$); namely,
the span of the first components of $R^J$ is equal to the Hopf subalgebra $H^J$ of
$D(H)^J,$ the span of the second components of $R^J$ is a Hopf
subalgebra of $D(H)^J$ which is naturally isomorphic to $(H^J)^*$ and the multiplication
map $H^J\ot (H^J)^*\to D(H)^J$ induces a Hopf algebra isomorphism between $D(H^J)$ and
$D(H)^J.$
We thus can identify $D(H)^J$ with $D(H^J),$ and using
this identification it is straightforward to check that the Drinfeld element of
$(D(H^J),R^J)$ is given by
\begin{equation}\label{uj}
u^J=Q^{-1}S(Q)u.
\end{equation}

We start with the following useful result.

\begin{proposition}\label{central}
Let $H$ be a finite-dimensional Hopf algebra, and let $g\in G(H)$ be such that
$gu^n$ is unipotent. Then $g=1.$
\end{proposition}

\begin{proof} First note that $g$ and $u$ commute by equation (\ref{ssq}).
Therefore,
$(1- gu^n)^N=0$ for some integer $N>0$ is equivalent to
$\sum_{k=0}^N (-1)^k \left(\begin{array}{c} N \\ k \end{array} \right)u^{nk}g^{k}=0,$
which in turn is equivalent to
$\sum_{k=0}^N (-1)^k \left(\begin{array}{c} N \\ k \end{array} \right)R_{nk}(g^{k}\ot
1)=0$ by Proposition \ref{eqdef}(3).
Now, apply $1\ot \varepsilon$ to the last equation to get $\sum_{k=0}^N (-1)^k
\left(\begin{array}{c} N \\ k \end{array} \right)g^{k}=0,$ i.e., \linebreak
$(1-g)^N=0.$ However, $1-g$ is semisimple, hence $g=1,$ as desired.
\end{proof}

As a corollary we can prove our first main result.
\begin{theorem}\label{tinv}
Let $H$ be a finite-dimensional Hopf algebra and let $J$ be a
twist for $H.$ Set $n:=\qexp(H).$ Then $u^n=(u^J)^n,$ and in particular
$\qexp(H^J)=n.$
\end{theorem}

\begin{proof}
Recall the formula
\begin{equation}\label{qsq}
\Delta(Q^{-1}S(Q))=J(Q^{-1}S(Q)\ot Q^{-1}S(Q))(S^2\ot S^2)(J^{-1})
\end{equation}
(see e.g. [AEGN]). Set $g:=S^{2n-1}(Q^{-1})S^{2n-2}(Q)\cdots
S^3(Q^{-1})S^2(Q)S(Q^{-1})Q.$
Since $S^{2n}=Id,$ we have by equation (\ref{qsq}) that $g$ is a grouplike element
in $H^J.$
By equation (\ref{uj}),
$g(u^J)^n=u^n,$ so $g(u^J)^n$ is unipotent, and the result follows from Proposition
\ref{central}.
\end{proof}

\begin{remark}
Theorem \ref{tinv} shows that $\qexp(H)$ is an invariant of the tensor category
of representations of $H.$ It would thus be interesting to develop a theory
of quasi-exponent for an arbitrary tensor category.
\end{remark}

Let us point out a direct corollary of Theorem \ref{tinv}, which may be of use in
future applications.
\begin{corollary}\label{ordgl}
Let $H$ be a finite-dimensional Hopf algebra and $J$ be a twist for $H.$ For any
grouplike
element $g$ in $H^J,$ $g^{\qexp(H)}=1.$
\end{corollary}

As another consequence we obtain the following.
\begin{corollary}\label{ddouble}
Let $H$ be a finite-dimensional Hopf algebra. Then $\qexp(D(H))=\qexp(H).$
\end{corollary}

\begin{proof}
By Proposition \ref{elem}(3), $\qexp(D(H))=\qexp(D(H)^*).$ Now, it is known that
$D(H)^*=(H^{op}\ot H^*)^R$ is obtained by twisting the Hopf algebra
$H^{op}\ot H^*$ using the $R-$matrix $R$ of $D(H).$ So by Theorem
\ref{tinv}, $\qexp(D(H))=\qexp(H^{op}\ot H^*).$ However, by Proposition
\ref{elem}(4),(7), $\qexp(H^{op}\ot H^*)=\qexp(H),$ and we are done.
\end{proof}

\end{section}

\begin{section}
{Quasi-exponent of finite-dimensional pointed Hopf algebras}
Recall that a Hopf algebra is called pointed if all its irreducible
corepresentations are $1-$dimensional (e.g. quantum universal enveloping algebras and
quantum groups at roots of unity). In this section we calculate the quasi-exponent of
finite-dimensional pointed Hopf algebras.

\begin{lemma}\label{aff}
Let $X$ be an affine connected algebraic variety, and let $H$ be a Hopf
algebra over the coordinate ring $\mathbb{C}[X]$ of $X,$ which is a finitely generated free
module of rank $r.$ For any $x\in X$ let
$I_x\subseteq \mathbb{C}[X]$ be the corresponding maximal ideal, and set $H_x:=H/I_xH;$ it
is a Hopf
algebra of dimension $r$ over $\mathbb{C}.$ Then $\qexp(H_x)$ is a constant which does
not depend on the point $x.$
\end{lemma}
\begin{proof}
Let $u_x$ be the Drinfeld element of $D(H_x).$ Then we know the eigenvalues of $u_x$
are roots of unity (see [EG]), hence the coefficients of the characteristic
polynomial of $u_x$ are algebraic numbers.
However, they must also continuously depend on $x.$ This
means they are constant, and so the spectrum of $u_x$ is independent on $x$, as desired.
\end{proof}

\begin{lemma}\label{filt}
Let $H$ be a $\mathbb{Z}_+ -$filtered Hopf algebra and let $\text{gr}H$ be its associated
graded Hopf
algebra. Then $\qexp(H)=\qexp(\text{gr}H).$
\end{lemma}
\begin{proof}
It is well known that there exists a Hopf algebra $\overline {H}$ over $\mathbb{C}[t]$
such that
$\overline {H}_0=\text{gr}H$ and $\overline {H}_t=H$ for all $t\ne 0.$ Thus the lemma
follows from Lemma \ref{aff}.
\end{proof}

\begin{proposition}\label{gr}
Let $H$ be a $\mathbb{Z}_+ -$graded Hopf algebra with zero part $H_0.$ Then
$\qexp(H)=l.c.m.(\qexp(H_0),|S^2|),$ where $|S^2|$ denotes the order of $S^2$ in $H.$
\end{proposition}
\begin{proof}
Set $n:=l.c.m.(\qexp(H_0),|S^2|),$ and let $R_0$ be the universal $R-$matrix
of $D(H_0).$ Since $\qexp(H_0)$ divides $n,$ one has
$\sum_{k=0}^N (-1)^k
\left(\begin{array}{c} N \\ k \end{array} \right)(R_0)_{nk}=0$ for some
positive integer $N.$ But $S^{2n}=Id$ on
$H_0,$ so this equation is equivalent to
$\sum_{k=0}^N (-1)^k
\left(\begin{array}{c} N \\ k \end{array} \right)(R_0)_{n}^k=0.$ That is, $(R_0)_{n}$ is
unipotent. This implies that $(R_0)_{n}$ is unipotent in the quotient algebra
$H_0/\Rad(H_0)\ot H_0^*/\Rad(H_0^*).$
However,
$$H/\Rad(H)\ot H^*/\Rad(H^*)=H_0/\Rad(H_0)\ot H_0^*/\Rad(H_0^*)$$
(since $H_{\ge
1},$ $H^*_{\ge 1}$ are nilpotents ideals), and moreover $R=R_0$ in this quotient algebra.
So $R_n$ is unipotent in $H/\Rad(H)\ot
H^*/\Rad(H^*),$ hence in $H\ot H^*,$ which is equivalent to saying that
$\sum_{k=0}^N (-1)^k
\left(\begin{array}{c} N \\ k \end{array} \right)R_{nk}=0$ for some
integer $N>0$ (since $S^{2n}=Id$ on $H$), hence
$n$ is divisible by $\qexp(H).$

Since by Proposition \ref{elem}, both $\qexp(H_0)$ and $|S^2|$ divide $\qexp(H),$ $n$
divides $\qexp(H),$ and we are done.
\end{proof}

\begin{theorem}\label{deg1}
Let $H$ be a finite-dimensional pointed Hopf algebra over $\mathbb{C}.$
Then $|S^2|$ divides $\exp(G(H)).$
\end{theorem}
\begin{proof}
Let $grH$ be the graded pointed Hopf algebra associated to $H$ with respect to
the coradical filtration of $H$ (see e.g. [M]). Then $|S_H^2|=|S_{grH}^2|$ and
$G(H)=G(grH).$ Therefore, it is sufficient to consider the graded case, so we will assume
that $H=\bigoplus H_i$ is graded.

We have a $\mathbb{Z}_+-$filtration of $H$ defined as follows: $F_i(H)$ is
defined to be the subalgebra in $H$ generated by $H_0,$ $H_1,\dots,$ $H_i.$
It is clearly a Hopf subalgebra in $H.$

Let $n:=\exp(G(H)),$ and let us show that $S^{2n}=Id$ on $F_i(H)$ by induction in $i.$

For $i=0,$ this is clear, since $F_0(H)=H_0=\mathbb{C}[G(H)].$
Suppose the statement is known for $i=k-1,$ and let us prove it for $i=k.$
Clearly, it is sufficient to show that $S^{2n}=Id$ on $H_k.$

Recall that $H_k$ is a bicomodule over $H_0.$ Indeed,
$\Delta:H_k\to \bigoplus_{i+j=k}H_i\ot H_j,$ so we can write $\Delta=\bigoplus_{i+j=k}
\Delta_{(i,j)},$ where $\Delta_{(i,j)}:H_k\to H_i\ot H_j.$ Now, the maps
$\Delta_{(0,k)},$ $\Delta_{(k,0)}$ define the $H_0-$bicomodule structure on $H_k.$
Therefore, since $H_0$ is cocommutative, we can restrict our attention to elements
of $H_k$
which belong to a $1-$dimensional bicomodule; i.e., elements $x\in H_k$
for which
$$\Delta(x)=g\otimes x+x\otimes h+\xi(x),$$
where $\xi(x)\in \bigoplus_{i=1}^{k-1}H_i\ot H_{k-i},$ and $g, h$ are grouplike
elements. In fact, we can assume $h=1$ (by multiplying $x$ by a grouplike element, if
necessary).

Now, apply $m(S\otimes Id)$ to this equation. Since $\varepsilon(x)=0,$
we get
$S(x)=-g^{-1}x$ modulo $F_{k-1}(H).$
Thus, $S^2(x)=g^{-1}xg$ modulo $F_{k-1}(H),$
and hence $S^{2n}(x)=x$ modulo $F_{k-1}(H)$ (as $g^n=1$).
Also by induction assumption, $S^{2n}=Id$ on $F_{k-1}(H).$
But $S$ is semisimple [R], hence $S^{2n}(x)=x,$ and we are done.
\end{proof}

As a corollary, we have the following weak version of a recent result by Radford
and Schneider [RS].
\begin{corollary}
Let $H$ be a finite-dimensional pointed Hopf algebra over $\mathbb{C}.$ Then $|S^2|$
divides $\dim(H).$
\end{corollary}
Radford and Schneider proved that $|S^2|$ divides $\dim(H)/|G(H)|$ [RS].

We can now prove our second main result.
\begin{theorem}\label{pointed}
Let $H$ be a finite-dimensional pointed Hopf algebra over $\mathbb{C}.$ Then
$\qexp(H)=\exp(G(H)).$
\end{theorem}
\begin{proof}
One has $\qexp(H)=\qexp(\text{gr}H)=l.c.m.(\exp(G(H)),|S^2|)=\exp(G(H))$ by Lemma
\ref{filt}, Proposition \ref{gr} and Theorem \ref{deg1}.
\end{proof}

As an application we can relate the exponents of the groups $G(H)$ and $G(H^J).$
\begin{corollary}\label{help3}
Let $H$ be a finite-dimensional pointed Hopf algebra over $\mathbb{C},$ and let
$J$ be a twist for $H.$
Then $\exp(G(H^J))$ divides $\exp(G(H)).$
\end{corollary}
\begin{corollary}\label{help4}
If two finite-dimensional pointed Hopf algebras $H_1,$ $H_2$ are twist equivalent,
then $\exp(G(H_1))=\exp(G(H_2)).$
\end{corollary}

\begin{example}
For example, Corollary \ref{help4} implies that two isocategorical groups (see [EG1]) have
the same exponent.
\end{example}

\begin{example}
Theorem \ref{pointed} can be applied to calculate the quasi-exponent
of quantum groups
at roots of unity (the finite-dimensional version). Indeed, if
${\mathfrak g}$ is a
finite-dimensional semisimple Lie algebra,
${\mathfrak b}_+$ its Borel subalgebra, and $q$ is a primitive
root of unity of order $l,$ then
$\qexp(U_q({\mathfrak b}_+))=\qexp(U_q({\mathfrak g}))=l.$
\end{example}

As a consequence we have the following result.

\begin{corollary}\label{qexpp}
For any twist $J$ of $U_q({\mathfrak g})$ and any grouplike element
$g$ of $U_q({\mathfrak g})^J,$ the order of $g$ divides $l.$
\end{corollary}

\begin{remark}\label{abelian}
We expect that moreover $G(U_q({\mathfrak g})^J)$ is abelian
(see Section 5).
\end{remark}

Let us conclude the section by calculating the quasi-exponent of
finite-dimensional
triangular Hopf algebras with the Chevalley property (see [AEG], [EG2]).

\begin{proposition}\label{chev}
Let $H$ be a finite-dimensional triangular Hopf algebra with the Chevalley
property over $\mathbb{C},$ and let $H_s:=H/\Rad(H)$ be its semisimple part. Then
$\qexp(H)=\exp(H_s).$
\end{proposition}
\begin{proof}
It is enough to consider the non-semisimple case. We have that $H_s^*$ is the zero
part in the coradical filtration of $H^*,$
therefore by Proposition \ref{gr},
$\qexp(H)=\qexp(H^*)=l.c.m.(\qexp(H_s^*),|S_H^2|).$ Now
by [AEG], $|S_H^2|=2.$
However, the Drinfeld element $u$ of $H$ is a grouplike element satisfying $u^2=1,$ $u\ne
1.$ Hence, $2=|S_H^2|$ divides $\exp(H_s),$ and we are done.
\end{proof}

\begin{remark}
Proposition \ref{chev} gives an efficient way of calculating $\qexp(H)$ since
by [EG3], $H_s=\mathbb{C}[G]^J$ for some finite group $G,$
hence $\exp(H_s)=\exp(G)$ by [EG, Theorem 3.3].
\end{remark}

\begin{remark}
The results of this paper, so far, extend without changes to the case when
the ground field
$\mathbb{C}$ is replaced by an algebraically closed field $k$ of odd
characteristic
relatively prime to the dimensions of the Hopf algebras involved.
More precisely, all
results remain true except Remark \ref{pch}(3).
(The assumption about odd characteristic is used when we
claim that $S^2$ is semisimple.)
\end{remark}

\end{section}

\begin{section}
{Applications to quantum groups at roots of unity}
Let ${\mathfrak g}$ be a finite-dimensional complex simple Lie algebra
and $p>2$ an odd number (not necessarily a prime). Assume also that $p$
is not divisible by $3$ if ${\mathfrak g}$ is of type $G_2.$
Let $q$ be a primitive $p-$th root of unity, and $U_q({\mathfrak g})$
denote the quantum group associated to ${\mathfrak g},$ generated by
$e_i,f_i,K_i,$ with the usual relations and the additional relations
$K_i^p=1,$ $e_i^p=f_i^p=0.$

The following theorem is essentially contained in [T] (see p. 404).

\begin{theorem}\label{tak} If $p$ is relatively prime to the
determinant of the Cartan matrix $A$ of ${\mathfrak g},$ then the Hopf
algebra $U_q({\mathfrak g})$ is simple, i.e. does not have non-trivial
Hopf ideals.
\end{theorem}

\begin{remark} We are grateful to N. Andruskiewitsch for explaining the
proof of this theorem to us.
\end{remark}

\begin{proof} Let $\phi: U_q({\mathfrak g})\to L$ be a non-trivial Hopf
algebra map of $U_q({\mathfrak g})$ into another finite-dimensional Hopf
algebra $L.$ Our job is to show that this map is injective.

First of all, let us show that $\phi(e_i)$ is non-zero for all $i.$ Indeed,
consider the set $I$ of all $i$ for which $\phi(e_i)=0.$
We claim that this set is empty. Indeed, if $i\in I$ then
$\phi([e_i,f_i])=0,$ so $\phi(K_i)=1$ (as $p$ is odd), hence
for any $j$ connected to $i,$ $\phi(e_j)=\phi(K_ie_jK_i^{-1})=
q^{a_{ij}}\phi(e_j),$ so $\phi(e_j)=0.$ So $I$ is either empty of everything.
But if $I$ is everything, we get $\phi(e_i)=0$ for all $i,$ hence
$\phi(K_i)=1$ for all $i,$ hence $\phi(f_i)=0$ for all $i,$
so $\phi$ is trivial. Contradiction.

Similarly, $\phi(f_i)$ is non-zero for all $i.$

Now, by Lemma 6.1.2 in [T], it is sufficient to check that
the map $\phi$ is injective in degree $1$ of the coradical
filtration of $U_q({\mathfrak g}).$ First let us check that it is so
in degree $0$ of this filtration, i.e. that $\phi: \mathbb{C}[G]\to L$
is injective, where $G$ is the group of grouplike elements generated by
$K_i.$ For this, it suffices to show that if $g\in G$ is non-trivial, then
so is $\phi(g).$ Assume the contrary, i.e. that $\phi(g)=1.$
Since $(p,\rm {det}(A))=1,$
there exists $i$ such that $ge_ig^{-1}=be_i,$ where $b\ne 1\in \mathbb{C}.$
Thus, $\phi(e_i)=0.$ Contradiction.

Now, the degree $1$ term of the coradical filtration has the form
$$\mathbb{C}[G]+\sum_i \mathbb{C}[G]e_i+\sum_i \mathbb{C}[G]f_i.$$
Assume that
$\phi$ is not injective in degree $1.$ Since $\rm {Ker}(\phi)$ is obviously
stable under conjugation by $G,$ there exists an element of $\rm
{Ker}(\phi)$ which is in $\mathbb{C}[G]e_i$ for some $i$ or $\mathbb{C}[G]f_i$
for some $i.$ Let us assume that the former is the case (the other case is
analogous). Then we have $\phi(\sum_{g\in G} a_g ge_i)=0.$ Calculating
the coproduct of this, one easily concludes that $\phi(e_i)=0.$
Contradiction. The theorem is proved.
\end{proof}

 From now on we suppose that $p$ is {\em prime.}

Let $Y$ be a non-trivial irreducible
representation of ${\mathfrak g},$ and assume that \linebreak
$p>\dim(Y).$

\begin{theorem}\label{thm1}
For any twist $J$ of $U_q({\mathfrak g}),$ the group
$G:=G(U_q({\mathfrak g})^J)$ is an elementary abelian $p-group.$
\end{theorem}

\begin{proof} By Corollary \ref{qexpp}, $\exp(G)=p.$ Thus, it remains
to show that $G$ is abelian. Consider the smallest abelian tensor
subcategory of $\Rep(U_q({\mathfrak g}))=\Rep(U_q({\mathfrak g})^J),$
which contains $Y.$ Since by Theorem \ref{tak}, $U_q({\mathfrak g})$ does
not have non-trivial Hopf algebra quotients, this category must coincide
with the whole $\Rep(U_q({\mathfrak g})).$ This means that any irreducible
representation occurs as a subquotient in the tensor
product of several copies of $Y.$

Now, since $G$ is a $p-$group, the dimensions of its irreducible
representations are powers of $p.$ This means that the restriction of
$Y$ to $G$ is a direct sum of $1-$dimensional representations of $G$
(as $p>\dim(Y)$).
Hence, the same holds for tensor powers of $Y,$ and therefore for
their subquotients; i.e., for all
irreducible representations of $U_q({\mathfrak g}).$
But then the regular representation of $U_q({\mathfrak g})$ is a sum of
$1-$dimensional representations for $G,$ which implies that $G$ is abelian.
\end{proof}

\begin{example} For ${\mathfrak g}:=sl(n),$ there is an $n-$dimensional
vector representation, so
the condition on $p$ is $p>n.$ For ${\mathfrak g}:=sp(2n)$ the condition is
$p>2n.$ For ${\mathfrak g}:=o(n),$ the condition is $p>n$ ($n\ge 7$).
\end{example}

\begin{conjecture}\label{con3}
Under the conditions of Theorem \ref{thm1}, the rank of
$G$ is at most the rank of the Lie algebra ${\mathfrak g}.$
\end{conjecture}

\begin{remark} There are examples (twist corresponding to 
Belavin-Drinfeld triples, see [EN]) where the rank of $G$ is actually less than
the rank of $g.$
\end{remark}

\begin{theorem}\label{thm4}
Conjecture \ref{con3} holds for classical the Lie algebras (series $A-D$).
\end{theorem}

\begin{proof}
Type $A_{n-1}=sl(n)$: Consider the restriction of the $n-$dimensional
irreducible representation $Y$ to $G.$ Then $Y=a_1+\cdots +a_n,$ where
$a_i$ are
characters of $G.$ Since $Y^{\otimes n}$ contains the trivial
representation as a subquotient
(quantum exterior power), there is at least one relation between
$a_1,\dots,a_n.$ This implies that the subgroup $K$ in $G^*$ generated by
$a_i$ is of rank $n-1.$ By the proof of Theorem \ref{thm1}, the characters of
$G$ which occur in any
irreducible representation of $U_q({\mathfrak g})$ must lie in $K,$ so
$K=G^*,$ and the rank of $G^*$ is at most $n-1,$ so we are done.

Types $C_n=sp(2n),$ $D_n=O(2n)$: Take the $2n-$dimensional irreducible
representation.
Since it is self-dual, the characters of $G$ that occur in it
are $a_i^{\pm 1},$ $i=1,\dots,n.$ But they generate $G^*,$ so the rank
of $G$ is at most $n.$

Type $B_n=O(2n+1):$ Take the $(2n+1)-$dimensional irreducible
representation. It is
self-dual, so the characters of $G$ that occur are: $1,a_i^{\pm 1},$
$i=1,\dots,n,$ so the rank of $G^*$ is at most $n.$

The theorem is proved.
\end{proof}

\end{section}

\begin{section}
{Acknowledgments}
P.E.\ was partially supported by the NSF grant DMS-9988796.
P.E.\ partially conducted his research for the
Clay Mathematics Institute as a Clay Mathematics Institute Prize Fellow.
S.G.\ thanks the Mathematics Department of MIT for the warm hospitality
during his visit. S.G.\ research was supported by the VPR - Fund at the
Technion and by the fund for the promotion of research at the Technion.
\end{section}

\bibliographystyle{ams-alpha}

\end{document}